\documentclass[11pt]{amsart}

\usepackage[english]{babel}
\usepackage{amsmath, amsfonts, amssymb, amsthm, epsfig, graphicx, url, hyperref, color, tikz, cancel,mathtools}
\usepackage{bbm}

\newtheorem{question}{Question}[section]

\newtheorem{theorem}{Theorem}[section]
\newtheorem{corollary}[theorem]{Corollary}
\newtheorem{lemma}[theorem]{Lemma}
\newtheorem{proposition}[theorem]{Proposition}

\newtheorem{remark}[theorem]{Remark}
\newcounter{theor}

\newtheorem{thm}[theor]{Theorem}

\mathtoolsset{showonlyrefs}

\usepackage[margin=0.9in]{geometry}

\def\s{\mathbb{S}}

\def\R{\mathbb{R}}

\newcommand{\dlat}{\mathrm{d}}

\def\esc#1{\left\langle #1\right\rangle}

\numberwithin{equation}{section}

\numberwithin{equation}{section}

\begin{document}
\title{Title}

\author[A. Malliaris]{Andreas Malliaris}
\address{Institut de Mathématiques de Toulouse (UMR 5219). Université de Toulouse
\& CNRS. UPS, F-31062 Toulouse Cedex 09, France}
\email{andreas.malliaris@math.univ-toulouse.fr}

\author[F. Marín Sola]{Francisco Marín Sola}
\address{Centro Universitario de la Defensa, Universidad Politécnica de Cartagena, C/ Coronel López Peña S/N, Base Aérea de San Javier, Santiago de La Ribera, 30720 Murcia, España. }
\email{francisco.marin7@um.es}

\thanks{ The first named author received support from the Graduate School EUR-MINT (State support managed by the National Research Agency for Future Investments program bearing the reference ANR-18-EURE-0023). The second named author is supported by the grant PID2021-124157NB-I00, funded by MCIN/AEI/10.13039/501100011033/``ERDF A way of making Europe'', as well as by the grant  ``Proyecto financiado por la CARM a través de la convocatoria de Ayudas a proyectos para el desarrollo de investigación científica y técnica por grupos competitivos, incluida en el Programa Regional de Fomento de la Investigación Científica y Técnica (Plan de Actuación 2022) de la Fundación Séneca-Agencia de Ciencia y Tecnología de la Región de Murcia, REF. 21899/PI/22''.}

\title{Quantitative improvements of functional inequalities under concavity properties}

\begin{abstract}
 A classical result of Hensley provides a sharp lower bound for the functional $\int_\R t^2f$, where $f$ is a non-negative, even log-concave function. In the context of studying the minimal slabs of the unit cube, Barthe and Koldobsky established a quantitative improvement of Hensley's bound. In this work, we complement their result in several directions. First, we prove the corresponding upper bound inequality for $s$-concave functions with $s\geq 0$. Second, we present a generalization of Barthe and Koldobsky's result for functionals of the form $\int_\R Nf\,\dlat\mu$, where $N$ is a convex, even function and $\mu$ belongs to a suitable class of positive Borel measures.  As a consequence of the employed methods, we obtain quantitative refinements of classical inequalities for $p$-norms and for the entropy of log-concave functions. Finally, we discuss both geometric consequences and probabilistic interpretations of our results. 
\end{abstract}
\maketitle

\section{Introduction}

Let $Q_n=\left[-\frac{1}{2},\frac{1}{2}\right]^n$ be the unit cube in $\mathbb{R}^n$. A classical result of Hadwiger \cite{hadwiger} states that among all central hyperplane sections of $Q_n$, i.e., $Q_n\cap\{x:\ \langle x,\theta\rangle=0\}$ for $\theta\in \s^{n-1}$, the ones parallel to some face have the minimal $(n-1)$-volume. It was thus natural to expect that for sufficiently small $t>0$, among all symmetric slabs of the cube, the canonical ones still have minimal volume. The fact that this $t$ would be independent of the dimension $n$ was conjectured by V. Milman and confirmed by Barthe and Koldobsky in \cite{BartheKoldobskySlabs} where, among other related results, it was shown that this conjecture holds at least when $t<\frac{3}{4}$. 

The approach of \cite{BartheKoldobskySlabs} was inspired by Hensley's proof of Hadwiger's result (see \cite{hensley}). It translates the geometrical problem into a functional inequality for decreasing log-concave functions on $[0,+\infty)$. Specifically, let $f: [0,+\infty)\to [0,+\infty)$ be an integrable,  decreasing log-concave function. It was shown in \cite{hensley} that 
\begin{equation}\label{e:hensley}
    f(0)^2\int_0^{+\infty}t^2 f(t)\ dt   \geq  \frac{1}{3}\left(\int_0^{+\infty} f(t)\ dt\right)^3.
\end{equation}
Barthe and Koldobsky provided the following quantitative improvement in \cite{BartheKoldobskySlabs}.
\begin{thm}[Barthe--Koldobsky] \label{t: Barthe Koldobsky functional}
    Let $f :[0,+\infty) \to [0,+\infty)$ be an integrable decreasing log-concave function. If $h>0$ is such that
    $$
    \int_0^h f \leq \frac{3}{4}\int_0^{+\infty} f,
    $$
    then
    \begin{equation}\label{e:Barthe_Koldobsky}
           3\int_0^{+\infty}t^2 f(t)\ dt \left(\int_0^h f(t)\ dt\right)^2  \geq h^2 \left(\int_0^{+\infty} f(t)\ dt\right)^3.
    \end{equation}
    Equality holds for $f(t)=c\chi_{_{[0,d]}}(t)$, with $c\geq 0$ and $d\geq h$.
\end{thm}

A classical result of Marshall, Olkim and Proschan \cite{MarshallOlkimProschan} (see also \cite[Proposition~5.3]{GNTconcentrationbook}) shows that the quantity 
$$
f(0)^2\int_0^{+\infty}t^2 f(t)\ dt
$$ 
admits an upper bound complementing \eqref{e:hensley}. This naturally leads to the question of whether a similar type of upper bound also exists for \eqref{e:Barthe_Koldobsky}. Our first result provides a positive answer in this direction.
\begin{theorem}\label{t:upper_bound intro}
    Let $f :[0,+\infty) \to [0,+\infty)$ be an integrable decreasing log-concave function. If $h>0$ is such that 
    $$
    \int_0^h f \leq (1-e^{-\sqrt{3}})\int_0^{+\infty} f,
    $$
    then
    \begin{equation}\label{e:Barthe_Koldobsky reverse}
           \int_0^{+\infty}t^2 f(t)dt \leq \frac{2h^2}{\left[-\log\left(1-\frac{\int_0^hf(t)dt}{\int_0^{+\infty}f(t)dt}\right)\right]^2} \int_0^{+\infty} f(t)dt.
    \end{equation}
     Equality holds for $f(t)=ce^{-\lambda t}\chi_{_{[0,+\infty)}}(t)$, with $c\geq 0$ and $\lambda>0$.
\end{theorem}

From a functional viewpoint, Theorem \ref{t:upper_bound intro} is the natural analogue of Theorem \ref{t: Barthe Koldobsky functional} for the upper bound. However, from a geometric point of view, this is not the case. More precisely, suppose that $K\subset\R^n$ is a symmetric convex body, that is convex, compact with non-empty interior, and symmetric with respect to the origin. As indicated above, such inequalities apply to the \emph{cross-section} function $t \mapsto |K\cap \{x\in \mathbb{R}^n:\ \langle x,\theta\rangle=t\}|_{n-1}$. According to Brunn's principle (see e.g. \cite[Theorem 1.2.2]{convexbook}) this function is $\frac{1}{n-1}$-concave, a condition stronger than log-concavity. Hence, since the extremal function in \eqref{e:Barthe_Koldobsky reverse} is merely log-concave, no sharp geometric result can be derived. To complement this, we establish the following functional inequality, which serves as a substitute for geometric considerations:
\begin{theorem}\label{theorem upper 1/n intro}
   Let $f:[0,+\infty)\to [0,+\infty)$ be a decreasing, integrable $\frac{1}{n-1}$-concave function and $n\in \mathbb{N}$. If $h>0$ is such that
\[
\int_0^h f \leq \left(1-\left(\frac{n}{n+2}\right)^{n}\right)\int_0^{+\infty}f,
\]
then
\begin{equation}\label{eq. theorem s-conc upper}
    \int_0^{+\infty}t^2f(t)\ dt\leq \frac{2h^2\int_0^{+\infty}f(t)dt}{ \left(n+1\right)\left(n+2\right)\left(1-\left(1-\frac{\int_0^hf(t)dt}{\int_0^{+\infty}f(t)dt}\right)^{\frac{1}{n}}\right)^2}.
\end{equation}
Equality holds for $f(t)=c\left(1-\lambda t\right)_+^{n-1}$, with $c\geq 0$ and $\lambda>0$.
\end{theorem}
We note that Theorem \ref{theorem upper 1/n intro} follows as a special case of a more general statement for $s$-concave functions, whose presentation is deferred to Section \ref{s:optimization}. Moreover, both Theorems \ref{t:upper_bound intro} and \ref{theorem upper 1/n intro} not only recover but also strengthen the corresponding classical inequalities in the limit $h \to 0^+$; see Remark \ref{remark improvement}.

A general family of inequalities for the case $h=0$ (that is, when $\int_0^h f/h$ is formally replaced by $f(0)$) was established by Fradelizi \cite{fradelizi}. There, sharp upper and lower bounds were derived for functionals of the form
\begin{equation} \label{Fra functional}
    \int_{\mathbb{R}} N(t) f(t)\,dt,
\end{equation}
where $f$ is log- or $s$-concave with barycenter at $0$ and $N$ is a non-negative, even convex function. This framework covers the class considered here, namely even $s$-concave functions on $\mathbb{R}$ with $s\ge 0$ and $N(t)=t^2$. Motivated by Theorem \ref{t: Barthe Koldobsky functional}, we also establish the following general result, which quantitatively improves all known lower bounds for \eqref{Fra functional}.
\begin{theorem} \label{main theorem intro}
     Let $N:[0,+\infty)\to [0,+\infty)$ be increasing and convex. Let $f:[0,+\infty)\to [0,+\infty)$ be a decreasing, integrable log-concave function. If $h>0$ is such that
     \begin{equation}
         \int_0^h f \leq \frac{1}{2}\int_0^{+\infty} f,
     \end{equation}
     then
     \begin{equation} \label{eq: statement of main theorem}
         \int_0^{+\infty} N(t)f(t)dt \geq  \frac{\int_0^h f}{h}\int_0^{\frac{h\int_0^{+\infty}f}{\int_0^hf}}N(t)dt.
     \end{equation}
     Equality holds for $f(t)=c\chi_{_{[0,d]}}(t)$, with $c\geq 0$ and $d\geq h$.
\end{theorem}
As aforementioned, \eqref{eq: statement of main theorem} improves the classical (lower bound) inequalities in \cite[Theorem 8]{fradelizi} when $\int_0^h f/h$ is replaced by $f(0)$, and it recovers them in the limit as $h\to 0^+$ (see Remark \ref{remark after main theorem}). Indeed, our contribution goes beyond this case; we establish a generalization of Theorem \ref{main theorem intro} for a broader class of measures supported on $[0,+\infty)$ thus obtaining a new family of inequalities as $h\to 0^+$. This is the content of Theorem \ref{main thm weighted}.

All the main results presented here follow from a general reduction scheme that was also used (for the lower bound and the Lebesgue case) in \cite{BartheKoldobskySlabs}. It should be compared with the powerful localization principle of Fradelizi and Guédon \cite{FradeliziGuedonGeneralLocaliz,FradeliziGuedonS-conc}, which is more general and has been widely used in many extremal problems concerning log-(or $s$-)concave  functions. Nevertheless, the reduction scheme that we employ can be easily adapted to other functionals. Specifically, by using this technique, we give new quantitative bounds for classical inequalities concerning the $p$-norms and the entropy of decreasing log-concave functions. This result reads as follows.
\begin{proposition}\label{p:p-norm}
    Let $f:[0,+\infty)\to [0,+\infty)$ be a decreasing, integrable log-concave function. For each $p\geq1$ there exists $\theta_p>0$, such that if $h>0$ satisfies
    \[
    \int_0^hf(t)dt\leq \theta_p \int_0^{+\infty} f(t)dt, 
    \]
    then
    \begin{equation} \label{p-norms estimate}
        \int_0^{+\infty} f(t)^pdt \leq \left( \frac{1}{h}\int_0^hf(t)dt\right)^{p-1} \int_0^{+\infty} f(t) dt.
    \end{equation}
    Equality holds for $f(t)=c\chi_{_{[0,d]}}(t)$, with $c\geq 0$ and $d\geq h$, and $\theta_p$ given by \eqref{theta_p in p-norms}.  
\end{proposition}

Other consequences of the main results are also explored, including probabilistic interpretations and geometric applications. From the probabilistic perspective, we obtain estimates for \emph{positive} median and the Laplace transform of an even log-concave random variable. In addition, we show that Theorem \ref{main theorem intro} yields the following quantitative refinement of Jensen’s inequality for random variables with log-concave tails.
\begin{corollary} \label{corollary jensen}
    Let $X$ be a non-negative random variable with log-concave tails. Let $N:[0,+\infty)\to[0,+\infty)$ with $N(0)=0$ and having convex and increasing derivative.
    If $h>0$ is such that
    $$
    \int_0^h\mathbb{P}(X\geq t)dt\leq \frac{1}{2}\mathbb{E}X,
    $$
    then
    \begin{equation}
        \mathbb{E}N(X)\geq \frac{\int_0^h\mathbb{P}(X\geq t)dt}{h} N\left(  \frac{h}{\int_0^h\mathbb{P}(X\geq t)dt} \mathbb{E}X \right).
    \end{equation}
\end{corollary}
In the geometric setting, we offer a new viewpoint on \cite[Theorem~1]{fradelizi} by working with \textit{slabs} of centrally symmetric convex bodies rather than central \textit{sections}. Combined with Theorem \ref{t: Barthe Koldobsky functional}, this leads to an interpretation of the recently solved slicing problem (see \cite{klartaglehec,guan,bizeul}) in terms of symmetric slabs instead of sections. Consequences for floating bodies are also discussed.

The paper is organized as follows. For the reader’s convenience, the proofs of the main theorems are divided into two sections: Section \ref{subsection reductions} reduces the problems to a simpler family of functions, while Section \ref{s:optimization} optimizes over that family. Applications of the developed techniques, together with probabilistic interpretations of the main results, are presented in Section \ref{section applications}. The aforementioned geometric applications are shown in Section \ref{section geometric}. Finally, some open problems and questions related with the main theorems are collected in Section \ref{section discussion}.

\section{Reductions}\label{subsection reductions} 
Here we prove that the treated extremal problems can be reduced to some simpler families of log-concave and $s$-concave functions, respectively. Let $\mu$ be a positive Borel measure on $[0,+\infty)$. Given $h>0$ and $V>u>0$, we define 
\begin{equation} \label{family of functions}
    \mathcal{F}_{h,u,V}(\mu):= \{ f:[0,+\infty)\to [0,+\infty)\text{, log-concave, decreasing,} \int_0^hfd\mu=u\text{ and }\int_0^{+\infty} f\ d\mu=V \},
\end{equation}
and
\begin{equation} \label{family of functions s-concave}
    \mathcal{F}_{h,u,V}^s(\mu):= \{ f:[0,+\infty)\to [0,+\infty)\text{, $s$-concave, decreasing,} \int_0^hfd\mu=u\text{ and }\int_0^{+\infty} f\ d\mu=V \},
\end{equation}
where $s>0$.

\begin{lemma} \label{lemma reduction upper}
     Let $N : [0,+\infty)\to [0,+\infty)$ be an increasing function. Let $h>0$ and $V>u>0$. Then,
   \[
        \sup\left\{\int_0^{+\infty}Nf\,\dlat\mu: f\in \mathcal{F}_{h,u,V}\right\}=
        \]\[\sup\left\{\int_0^{+\infty}Nf\,\dlat\mu: f\in \mathcal{F}_{h,u,V} \text{ and } f(t) = c\chi_{_{[0,d]}}(t) + c e^{-\lambda(t-d)}\chi_{_{[d,+\infty)}}(t), d\leq h \text{ and } c,\lambda \geq 0\right\}.
   \]
\end{lemma}
\begin{proof}
    Let $f=\exp(\phi)\in \mathcal{F}_{h,u,V},$ where $\phi$ is concave and decreasing on $[0,+\infty)$. We denote below by $\phi'_-(s)$ and $\phi'_+(s)$ the left and right derivative respectively at the point $s$. For each $p\geq -\phi_+'(h)$ and $d\leq h$ consider the function 
    \[g_{p,d}(t)=f(h)e^{-p(d-h)}\chi_{_{[0,d]}}(t)+f(h)e^{-p(t-h)}\chi_{_{[d,+\infty)}}(t),\]
    where $t\geq 0$. Note that the function $g_{p,d}$ is of the desired form. For each $t\geq h$ and $d\leq h$ we have that $g_{-\phi_+'(h),d}(t)\geq f(t)$, by concavity, while trivially $0=g_{+\infty,d}(t)\leq f(t)$. This implies the existence of $p\in [-\phi_+'(h),+\infty]$ such that 
    \[
    \int_h^{+\infty}f(t)d\mu(t)= \int_h^{+\infty}g_{p,d}(t)d\mu(t).
    \]
    Moreover, concavity implies that unless the two functions coincide there is a unique intersection point, say $t_+$, for which $g_{p,d}\leq f$ on $[h,t_+]$ and $g_{p,d}\geq f$ on $[t_+,+\infty)$. This was done independently of the choice of $d$. For $t\leq h$, since $\phi'_-(h)\geq \phi'_+(h)\geq -p$, $g_{p,0}(t)\geq f(t)$ by concavity. On the other hand, since $f$ is decreasing, $f\geq g_{p,h}$ on $[0,h]$. By continuity, there exists $d\leq h$ such that 
    \[
    \int_0^hf(t) \;d\mu(t)=\int_0^h g_{p,d}(t)\;d\mu(t).
    \]
    Concavity again ensures that unless the functions coincide there exists a unique $t_-\in [0,h]$ such that $f\geq g_{p,d}$ on $[0,t_-]$ and $f\leq g_{p,d}$ on $[t_-,h]$. Finally, from the above discussion and the fact that $N$ is increasing we get
    \[
    \int_0^{+\infty}N(t)(f(t)-g_{p,d}(t))d\mu(t)=
    \]
    \[
    \int_0^{t_-}N(t)(f(t)-g_{p,d}(t))d\mu(t) +\int_{t_-}^hN(t)(f(t)-g_{p,d}(t))d\mu(t)\]\[ + \int_{h}^{t_+}N(t)(f(t)-g_{p,d}(t))d\mu(t) +\int_{t_+}^{+\infty}N(t)(f(t)-g_{p,d}(t))d\mu(t)
    \]
    \[
    \leq N(t_-)\int_0^{t_-}(f(t)-g_{p,d}(t))d\mu(t)+N(t_-)\int_{t_-}^h(f(t)-g_{p,d}(t))d\mu(t)
    \]
    \[ +N(t_+) \int_{h}^{t_+}(f(t)-g_{p,d}(t))d\mu(t) +N(t_+)\int_{t_+}^{+\infty}(f(t)-g_{p,d}(t))d\mu(t)=0.
    \]
    \end{proof}

\begin{lemma} \label{lemma reduction lower}
     Let $N:[0,+\infty)\to [0,+\infty)$ be an increasing function. Let $h>0$ and $V>u>0$. Then,
\[
\inf\left\{ \int_0^{+\infty} Nfd\mu :  f\in \mathcal{F}_{h,u,V}(\mu) \right\} =\]\[ \inf\left\{ \int_0^{+\infty} Nf\ d\mu :  f\in \mathcal{F}_{h,u,V}(\mu)   \text{ such that } f(t)=c \chi_{_{[0,d]}}e^{-at}, \text{  for  } a,c,d\geq 0, d\geq h\right\}.
\]
\end{lemma}
\begin{proof}
    Let $f=\exp(\phi)\in \mathcal{F}_{h,u,V},$ with $\phi$ concave and decreasing on $[0,+\infty)$. As before, we use $\phi'_-(s)$ and $\phi'_+(s)$ for the left and right derivative at the point $s$, respectively. For each $w\in [\phi'_-(h),0]$ and $v\in [h,+\infty]$ we consider the following concave and decreasing function on $[0,+\infty)$:
\[
\phi_{w,v}=\begin{cases}
w(t-h)+\phi(h) & t\leq v, \\
 -\infty & t>v. \\
\end{cases}
\]
It is clear that $\phi_{w,v}$ is of the desired form. First, we will find suitable values for $w$ and $v$ so that $e^{\phi_{w,v}}\in \mathcal{F}_{h,u,V}(\mu)$ and 
\begin{equation}\label{e:Lema_reduction_inf_1}
    \int_0^{+\infty} Ne^\phi d\mu \geq \int_0^{+\infty} Ne^{\phi_{w,v}} d\mu.
\end{equation}
For each $t\in [0,h]$ (note that $v\geq h$ so on that interval $\phi_{w,v}$ is independent of $v$) we have that $\phi_{0,v}=\phi(h)\leq \phi(t)$, since $\phi$ is decreasing and $\phi_{\phi'_-(h),v}(t)=\phi'_-(h)(t-h)+\phi(h)\geq \phi(t)$ by concavity. Hence, by continuity, we can find $\lambda \in [\phi'_-(h),0]$ such that (for each such $v$)
\[
\int_0^he^{\phi_{\lambda,v}(t)}d\mu(t)=\int_0^he^{\phi(t)}d\mu(t).
\]
Moreover, unless $\phi$ and $\phi_{\lambda,v}$ coincide on $[0,h]$, there exists a unique point $t_1\in (0,h)$ such that $\phi_{\lambda,v}\geq \phi$ on $[0,t_1]$ and $\phi_{\lambda,v}\leq \phi$ on $[t_1,h]$. We fix this $\lambda$.

Next, for each $t\in [h,+\infty)$ we obviously have $\phi_{\lambda,h}(t)\leq \phi(t)$, and $\phi_{\lambda,+\infty}(t)=\lambda(t-h)+\phi(h)\geq \phi(t)$ again by concavity, since $\lambda\geq \phi'_-(h)\geq\phi'_+(h)$. Hence, by continuity, there exists $d>h$, such that 
\[
\int_h^{+\infty} e^{\phi_{\lambda,d}(t)}d\mu(t)=\int_h^{+\infty} e^{\phi(t)}d\mu(t).
\]
Fixing also $d$, this makes sure that $e^{\phi_{\lambda,d}}\in \mathcal{F}_{h,u,V}(\mu)$. Now, unless $\phi_{\lambda,d}$ and $\phi$ coincide on $[h,+\infty)$, we have that $\phi_{\lambda,d}\geq \phi$ on $[h,d]$ (by concavity) and $\phi_{\lambda,d}\leq \phi$ on $[d,+\infty)$.

Lastly, combining all the above with the fact that $N$ is increasing, one can confirm \eqref{e:Lema_reduction_inf_1} following the strategy of Lemma \ref{lemma reduction upper}.
\end{proof}

\begin{figure}[h]
\centering
\begin{tikzpicture}[scale=1.5]
\draw[very thin] (-0.2,0) -- (3.6,0);
\draw[very thin] (0,3) -- (0,-0.2);
\draw[semithick,blue]
(0, 2.2) to [out=0, in=100] (3,0);
\draw[semithick, red,  domain=0:1.62]
plot({\x},{2});
\draw[semithick, red, domain=1.62:3]
plot({\x},{1.19*(-\x+ 1.62)+2});
\draw (0,0) circle (0.25pt) (-0.2,-0.06) node[below]{$0$};
\draw[very thin, dashed]
(1.89,1.68)--(1.89,0);
\draw[fill] (1.89,1.68) circle (0.9pt);
\draw (1.89,0) circle (0.25pt) (1.89,-0.02) node[below]{$h$};
\draw (3,0) circle (0.25pt) (3,-0.02) node[below]{$b$};
\draw  (0.57,1.5)  node{$\textcolor{red}{\tilde{g}_{p,d}^s}$};
\draw  (2.5,1.8)node{$\textcolor{blue}{f^{s}}$};
\end{tikzpicture}
\caption{The functions $f^s$ and $\tilde{g}_{p,d}^s$ in Lemma \ref{lemma reduction upper s}.}
\label{f:f^s_g^s}
\end{figure}
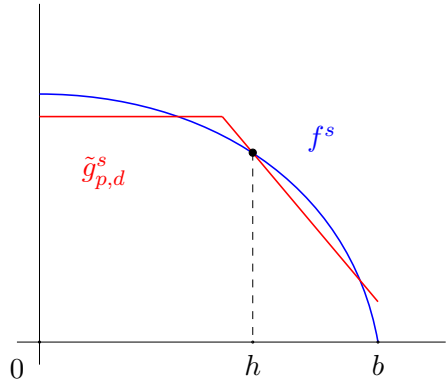

We will need a final reduction lemma, which is the analog of Lemma \ref{lemma reduction upper} for $s$-concave functions.
\begin{lemma}\label{lemma reduction upper s}
     Let $N : [0,+\infty) \to [0,+\infty)$ be an increasing function. Let $h>0$ and $V>u>0$. Then, 
   \[
        \sup\left\{\int_0^{+\infty}Nf\,\dlat\mu: f\in \mathcal{F}^s_{h,u,V}\right\}=\]\[
        \sup\left\{\int_0^{+\infty}Nf\,\dlat\mu: f(t) = c\chi_{_{[0,d]}}(t) + c \left(\frac{b-t}{b-d}\right)^{1/s}\chi_{_{[d,b]}}(t)\in \mathcal{F}^s_{h,u,V},\text{ with }d\leq h\right\}.
   \]
\end{lemma}
\begin{proof}
      Let $f\in \mathcal{F}^s_{h,u,V}$ and suppose that $supp(f) = [0,b]$. Consider the $s$-concave function $g_p(t) = f(h)\bigl(1 + p(t-h)\bigr)^{1/s}$ on $[h,b]$. There exists $p_0$ such that $\int_h^{b}g_{p_0}(t)\,\dlat \mu(t) = \int_h^{b}f(t)\,\dlat \mu (t)$. Since $g_{p_0}(h) = f(h)$ and by monotonicity, there exists a unique crossing point $t_1 \in (h,b)$, unless the functions coincide (see Figure \ref{f:f^s_g^s}). Moreover, using that $N$ is increasing, for such $p_0$ we have that
    \begin{align*}
        \int_h^{b} N(t)\bigl(g_{p_0}(t) - f(t)\bigr) \,\dlat \mu(t) =  \int_h^{b} \bigl(N(t) - N(t_1)\bigr)\bigl(g_{p_0}(t) - f(t)\bigr) \,\dlat \mu(t) \geq 0.
    \end{align*}
   For each $d\leq h$ consider the function 
   $$
   \tilde{g}_{p_0,d}(t) = f(h)\bigl(1 + p_0(d-h)\bigr)^{1/s}\chi_{_{[0,d]}}(t) + f(h)\bigl(1 + p_0(t-h)\bigr)^{1/s}\chi_{_{[d,b)}}(t).
   $$
   Again, there exists by continuity $d_0 \leq h$ such that $\int_0^{h}\tilde{g}_{p_0,d_0}(t)\,\dlat\mu(t) = \int_0^{h}f(t)\,\dlat\mu(t)$. Again, there exists a unique crossing point $t_0\in(0,h)$ unless the two functions coincide and 
    \begin{align*}
        \int_0^{h}N(t)\bigl(\tilde{g}_{p_0,d_0}(t) - f(t)\bigr)\,\dlat\mu(t) = \int_0^{h} \bigl(N(t) - N(t_0)\bigr)\bigl(\tilde{g}_{p_0,d_0}(t) - f(t)\bigr)\,\dlat\mu(t) \geq 0,
    \end{align*}
    which gives the claim.
\end{proof}

\section{Optimization}\label{s:optimization}
This section is devoted to the optimization part of the main theorems' proof. We start by that of Theorems \ref{t:upper_bound intro} and \ref{theorem upper 1/n intro}. Although the corresponding reduction lemmas are given for a general positive Borel measure and a increasing function $N :[0,+\infty) \to [0,+\infty)$, note that for the upper bounds we only use the cases of the Lebesgue measure and $N(t)=t^2$. 
 
\begin{proof}[Proof of Theorem \ref{t:upper_bound intro}]
    
By Lemma \ref{lemma reduction upper} we know that
\[
\sup \left\{ \int_0^{+\infty}t^2f(t)dt : f\in \mathcal{F}_{h,u,V}  \right\}= \]\[\sup\{c\int_0^dt^2dt+c\int_d^{+\infty}t^2e^{-\lambda(t-d)}dt :  \  c\chi_{_{[0,d]}}(t) + c e^{-\lambda(t-d)}\chi_{_{[d,+\infty)}}(t)\in  \mathcal{F}_{h,u,V}, d\leq h \}.
\]
So let $f(t)=c\chi_{_{[0,d]}}(t) + c e^{-\lambda(t-d)}\chi_{_{[d,+\infty)}}(t)\in  \mathcal{F}_{h,u,V}$, with $d\leq h$. Then
\[
u=\int_0^hf(t)dt\implies u=cd+\frac{c}{\lambda}-\frac{c}{\lambda}e^{-\lambda(h-d)},
\]
\[
V=\int_0^{+\infty}f(t)dt\implies V=cd+\frac{c}{\lambda}.
\]
Setting $x:=d\lambda$ and also denoting $\Delta:= 1-\frac{u}{V}$, we get:
\[
\frac{c}{\lambda}=\frac{V}{x+1},
\]
and
\[
V-u=\frac{c}{\lambda}e^{-\lambda h+x}\iff \lambda = \frac{x-\log[\Delta(x+1)]}{h},
\]
while 
\[
d=\frac{hx}{x-\log[\Delta(x+1)]}.
\]
Now, calculating the second moment we have 
\begin{align*}
    \int_0^{+\infty}t^2f(t)dt&=\frac{1}{3}cd^3+\frac{cd^2}{\lambda}+\frac{2V}{\lambda^2}\\
    &=Vh^2\frac{x^3+3x^2+6x+6}{3(x+1)[x-\log[\Delta(x+1)]]^2}.
\end{align*}
Observe that for $x=0$, the quantity becomes $2Vh^2/[-\log\Delta]^2$, which is the desired one. After simplifying, it suffices to prove that
 \[
 \frac{x^3+3x^2+6x+6}{3(x+1)[x-\log[\Delta(x+1)]]^2}\leq \frac{2}{[-\log\Delta]^2}.
 \]
 Or, after setting $\alpha^{-1}:= \log\frac{1}{\Delta}$, 
 \[
 (6x+6)[\alpha(x-\log(x+1))+1]^2-x^3-3x^2-6x-6\geq 0.
 \]
 Let $H(x):=(6x+6)[\alpha^2(x-\log(x+1))^2+2\alpha(x-\log(x+1))]-x^3-3x^2$. Since $H(0)=0$, it is enough to prove that $H$ is increasing. We calculate \[H'(x)=(18\alpha^2-3)x^2+(24\alpha-6)x-24\alpha^2x\log(x+1)+6\alpha^2\log(x+1)^2-12\alpha\log(x+1).\] Observe that $H'(0)=0$. Hence, we will show that $H'$ is increasing. Equivalently, we will prove that
 $$
 H''(x)=(36\alpha^2-6)x+(24\alpha-6)-24\alpha^2\log(x+1)-24\alpha^2\frac{x}{x+1}+12\alpha^2\frac{\log(x+1)}{x+1}-\frac{12\alpha}{x+1}\geq 0
 $$ 
or rearranging,  
 \begin{equation}\label{e:Th1.1}
 (6\alpha^2-1)x^2+(2\alpha^2+4\alpha-2)x+2\alpha-1-4\alpha^2(x+1)\log(x+1)+2\alpha^2\log(x+1)\geq 0.
 \end{equation}
At $x=0$ the above expression equals $2\alpha-1$, and so this is where a first restriction for $\alpha$, i.e., $\alpha \geq 1/2$ comes. Yet again, we will verify that \eqref{e:Th1.1} is an increasing function on $x$. Taking the derivative,  we would have to confirm that 
\begin{equation}\label{e:Th1.1_2}
    (12\alpha^2-2)x-2(\alpha-1)^2-4\alpha^2\log(x+1)+\frac{2\alpha^2}{x+1}\geq 0.
\end{equation}
Note that by the above restriction $\alpha \geq 1/2$, the value at $x=0$ of \eqref{e:Th1.1_2} is non-negative. Thus, our aim is once more reduced to check that this function is increasing on $x$. Taking derivatives this is equivalent to 
\[
(12\alpha^2-2)-4\alpha^2\frac{1}{x+1}-2\alpha^2\frac{1}{(x+1)^2}\geq 0.
\]
In order to ensure this, we will impose that $\alpha\geq \sqrt{3}^{-1}$, which is where the final restriction for $\alpha$ comes. Thus, taking $\frac{u}{V}\leq 1-e^{-\sqrt{3}}$ and rolling back, the claim follows.
\end{proof}

\begin{remark} \label{remark improvement}
As $h\to 0^+$, the right handside of \eqref{e:Barthe_Koldobsky reverse} converges to 
\[
\frac{2}{f(0)^2}\left(\int_0^{+\infty}f(t)dt\right)^3,
\]
and thus it recovers the corresponding upper bound from Hensley's inequality. Moreover writing
   \begin{equation}
       \frac{h}{-\log\left(1-\frac{\int_0^hf(t)dt}{\int_0^{+\infty}f(t)dt}\right)}= \frac{h}{-\log\left(\int_h^{+\infty}f(t)dt\right)+\log\left(\int_0^{+\infty}f(t)dt\right)},
   \end{equation}
   we observe that this expression is decreasing in $h$. Indeed, since $f$ is log-concave, its tails are log-concave as well and  the function $h\mapsto -\log\left(\int_h^{+\infty}f(t)dt\right)$ is convex. Therefore, for any $h>0$
   \[
    \frac{h}{-\log\left(1-\frac{\int_0^hf(t)dt}{\int_0^{+\infty}f(t)dt}\right)}\leq \frac{\int_0^{+\infty}f(t) dt}{f(0)},
   \]
   which means that \eqref{e:Barthe_Koldobsky reverse} is an improvement of the classical bound given by Hensley.
   
\end{remark}

For the case of $s$-concave functions the following holds recovering Theorem \ref{theorem upper 1/n intro} for $s=\frac{1}{n-1}$.
\begin{theorem}\label{theorem s-conc square}
    Let $f:[0,+\infty)\to [0,+\infty)$ be a decreasing, integrable $s$-concave function with $s\geq 0$. If $h>0$ is such that 
\[
\int_0^h f \leq \left(1-\left(\frac{s+1}{3s+1}\right)^{\frac{1}{s} + 1} \right)\int_0^{+\infty}f,
\]
then
\begin{equation}\label{eq. theorem s-conc upper}
    \int_0^{+\infty}t^2f(t)\ dt\leq \frac{2h^2\int_0^{+\infty}f(t)dt}{ \left(\frac{1}{s}+2\right)\left(\frac{1}{s}+3\right)\left(1-\left(1-\frac{\int_0^hf(t)dt}{\int_0^{+\infty}f(t)dt}\right)^{\frac{s}{s+1}}\right)^2}.
\end{equation}
Equality holds for $f(t)=c\left(1-\lambda t\right)_+^{\frac 1s}$, with $c\geq 0$ and $\lambda>0$.
\end{theorem}
Since the proof is long and tedious, it is deferred to the Appendix. Observe that as $s\to 0^+$ the right-hand side of \eqref{eq. theorem s-conc upper} converges to the right-hand side of \eqref{e:Barthe_Koldobsky reverse}. Moreover, arguing as in Remark \ref{remark improvement} we can see that \eqref{eq. theorem s-conc upper} improves upon the classical case, i.e., that of $h\to 0^+$.

\begin{remark}
    Both theorems of this section seem to generalize, via the same method, to (at least integer) moments $\int_0^{+\infty}t^pf(t)\ dt$. In lack of a specific application or a general principle, we decide not to delve in these very technical proofs. 
\end{remark}

\subsection{Lower bound optimization}
 Let $p\in [0,1)$ and $\lambda\geq 0$. We consider the measure on $[0,+\infty)$ given by $d\mu_{p,\lambda}(t)=t^{-p}e^{-\lambda t}dt$, and $\Phi_{p,\lambda}(t)=\mu_{p,\lambda}\left([0,t]\right)$ its distribution function. Let 
\[
\mathcal{I}:=\inf\left\{\int_0^{+\infty}N(t)f(t)\ d\mu_{p,\lambda}(t): f\in \mathcal{F}_{h,u,V}(\mu_{p,\lambda}) \right\},
\]
where we will denote below, unless otherwise stated, $V=\int_0^{+\infty}f\ d\mu_{p,\lambda}$ and $u=\int_0^{h}f\ d\mu_{p,\lambda}$.

In contrast with previous results for the corresponding upper bound inequality, for the lower bound we are able to optimize for the above class of measures and any convex increasing $N$. Our main result, which for $p=\lambda=0$ gives Theorem \ref{main theorem intro}, reads as follows:
\begin{theorem} \label{main thm weighted}
Let $N:[0,+\infty)\to [0,+\infty)$ be increasing and convex. Let $f:[0,+\infty)\to [0,+\infty)$ be a decreasing, integrable log-concave function. If $h>0$ is such that 
\[
\int_0^h f \ d\mu_{p,\lambda} \leq \frac{1}{2}\int_0^{+\infty}f\ d\mu_{p,\lambda},
\]
then 
    \begin{equation} \label{equation main theorem}
        \int_0^{+\infty} N(t)f(t)d\mu_{p,\lambda}(t)\geq \frac{u}{\Phi_{p,\lambda}(h)} \int_0^{\Phi_{p,\lambda}^{-1}\left(\frac{V}{u}\Phi_{p,\lambda}(h)\right)}N(t)d\mu_{p,\lambda}(t).
    \end{equation}
    Equality holds for $f(t)=c\chi_{_{[0,d]}}(t)$.
\end{theorem}

We will need the following auxiliary lemma:
\begin{lemma} \label{lemma properties of function}
    Let $F_p(x)=\int_0^{x}t^{-p}e^{-t}dt, x\geq 0$, and set $G_p:=F_p^{-1}$. Then, 
    \begin{enumerate}
        \item $G_p$ is increasing,
        \item $G_p(t)^{-p}e^{-G_p(t)}$ is decreasing,
        \item $G_p(t)^{-p+1}e^{-G_p(t)}$ is concave.
    \end{enumerate}
\end{lemma}
\begin{proof}
The first claim is obvious. The last two follow from the fact that $\int_0^{G_p(x)}t^{-p}e^{-t}dt=x$ and thus $G_p'(x)G_p(x)^{-p}e^{-G_p(x)}=1$. Computing the derivatives of the two functions along with the first item gives the rest assertions.
\end{proof}

\begin{lemma} \label{lemma optimization}
In the notation of Theorem \ref{main thm weighted} we have:
    \begin{equation}
        \mathcal{I}= \inf\left\{\int_0^V N\left(h \frac{G_p(\frac{tx}{V})}{G_p(\frac{ux}{V})}  \right)dt :    x\in \left(\frac{V}{u}F_p(\lambda h),\Gamma(-p+1)\right)  \right\}
    \end{equation}
\end{lemma}
\begin{proof}
By Lemma \ref{lemma reduction lower}, we know that the infimum $\mathcal{I}$ is achieved for functions of the form $f(t)=c \chi_{_{[0,d]}}e^{-at}$. For $z:=a+\lambda$ $(\geq \lambda)$, we set 
\[s= \int_0^{t}f\;d\mu_{p,\lambda}= \frac{c}{z^{-p+1}}F_p(zt)\iff t= \frac{1}{z}G_p\left( \frac{z^{-p+1}}{c}s \right)\]
and so 
    \[
    \int_0^{+\infty}N(t)f(t)d\mu_{p,\lambda}(t) = \int_0^{+\infty}N\left( \frac{1}{z}G_p\left( \frac{z^{-p+1}}{c}s \right) \right)ds.
    \]
    We have also that \[
    u=c\int_0^h e^{-zt}t^{-p}dt \quad \mathrm{and} \quad V=c\int_0^d e^{-zt}t^{-p}dt.
    \] 
    Hence,
    \[
    u=\frac{c}{z^{-p+1}}F_p(zh) \quad \text{and} \quad V= \frac{c}{z^{-p+1}}F_p(zd).
    \]
    Now, setting $x= \frac{V}{u}F_p(zh)\geq \frac{V}{u}F_p(\lambda h)$, we have that $x=F_p(zd)$ and thus $x< \Gamma(-p+1)$. Rewriting the parameters involved in the above expression as 
    \[
    z=\frac{1}{h} G_p\left({\frac{u}{V}x}\right) \quad \text{and} \quad \frac{z^{-p+1}}{c}=\frac{x}{V},
    \]
    we deduce the desired claim.
\end{proof}

\begin{proof}[Proof of Theorem \ref{main thm weighted}]
We may assume after normalizing that $V=1$. From Lemma \ref{lemma optimization}, we define 
\[
K_u(x)=\int_0^1 N\left(h \frac{G_p(tx)}{G_p(ux)}  \right)dt,\;\text{for } x \in \left(\frac{V}{u}F_p(\lambda h),\Gamma(-p+1)\right)
\]
and  $u\in (0,\frac{1}{2}]$ fixed. We will assume that $N$ is differentiable, which also implies continuity of its derivative. It is thus enough to prove that $K'_u(x)\geq 0$ for each $x$. Equivalently, after differentiating and removing an extra term of $G_p(ux)$ and $h$, we would like to show that for each $x$
    \[
    \int_0^1 N'\left(h \frac{G_p(tx)}{G_p(ux)}  \right)\left[tG_p'(tx)-\frac{uG_p(tx)G_p'(ux)}{G_p(ux)}\right]dt \geq 0.
    \]
    Setting $tx=s$ and $ux=r\leq x/2$, we equivalently want (after multiplying by $x^2$)
    \[
    H(x):=\int_0^x N'\left(h \frac{G_p(s)}{G_p(r)}  \right)\left[sG_p'(s)-\frac{rG_p(s)G_p'(r)}{G_p(r)}\right]dt \geq 0
    \]
    for all $x \in \left(\frac{V}{u}F_p(\lambda h),\Gamma(-p+1)\right)$. The computations in Lemma \ref{lemma properties of function} show that 
\[
H(x) = \int_0^x N'\left(h \frac{G_p(s)}{G_p(r)}  \right)\left[s\frac{e^{G_p(s)}}{G_p(s)^{-p}}-\frac{rG_p(s)e^{G_p(r)}}{G_p(r)^{-p+1}}\right]dt .
\]
Observe that by the concavity of $G_p(t)^{-p+1}e^{-G_p(t)}$,
\[
s\frac{e^{G_p(s)}}{G_p(s)^{-p}}-\frac{rG_p(s)e^{G_p(r)}}{G_p(r)^{-p+1}}\geq 0 \iff s\geq r.
\]
Since $r\leq x/2$, we may rewrite $H(x)$ as follows:

\[
H(x)=\int_0^{r}N'\left(h \frac{G_p(s)}{G_p(r)}  \right)\left[s\frac{e^{G_p(s)}}{G_p(s)^{-p}}-\frac{rG_p(s)e^{G_p(r)}}{G_p(r)^{-p+1}}\right]dt\]
\[+\int_{r}^{2r}N'\left(h \frac{G_p(s)}{G_p(r)}  \right)\left[s\frac{e^{G_p(s)}}{G_p(s)^{-p}}-\frac{rG_p(s)e^{G_p(r)}}{G_p(r)^{-p+1}}\right]dt
\]
\[
+\int_{2r}^xN'\left(h \frac{G_p(s)}{G_p(r)}  \right)\left[s\frac{e^{G_p(s)}}{G_p(s)^{-p}}-\frac{rG_p(s)e^{G_p(r)}}{G_p(r)^{-p+1}}\right]dt.
\]
The first term is negative while the last two are positive by the previous observation. It will be sufficient thus to show that
\[
\int_{r}^{2r}N'\left(h \frac{G_p(s)}{G_p(r)}  \right)\left[s\frac{e^{G_p(s)}}{G_p(s)^{-p}}-\frac{rG_p(s)e^{G_p(r)}}{G_p(r)^{-p+1}}\right]dt\geq \int_0^{r}N'\left(h \frac{G_p(s)}{G_p(r)}  \right)\left[-s\frac{e^{G_p(s)}}{G_p(s)^{-p}}+\frac{rG_p(s)e^{G_p(r)}}{G_p(r)^{-p+1}}\right]dt.
\]

Setting $2r-t=s$ in the first integral (and then renaming the variable) this is equivalent to
\[
\int_{0}^{r}N'\left(h \frac{G_p(2r-s)}{G_p(r)}  \right)\left[(2r-s)\frac{e^{G_p(2r-s)}}{G_p(2r-s)^{-p}}-\frac{rG_p(2r-s)e^{G_p(r)}}{G_p(r)^{-p+1}}\right]dt\geq\]\[ \int_0^{r}N'\left(h \frac{G_p(s)}{G_p(r)}  \right)\left[-s\frac{e^{G_p(s)}}{G_p(s)^{-p}}+\frac{rG_p(s)e^{G_p(r)}}{G_p(r)^{-p+1}}\right]dt.
\]
This is true pointwise. Indeed, since $2r-s\geq s$ ($s\leq r$), $N$ is convex and $G_p$ increasing, we have that
$$
N'\left(h \frac{G_p(2r-s)}{G_p(r)}\right)\geq N'\left(h \frac{G_p(s)}{G_p(r)}  \right).
$$
It remains now to show that for each $s\leq r$ it holds that
\[
(2r-s)\frac{e^{G_p(2r-s)}}{G_p(2r-s)^{-p}}-\frac{rG_p(2r-s)e^{G_p(r)}}{G_p(r)^{-p+1}} \geq -s\frac{e^{G_p(s)}}{G_p(s)^{-p}}+\frac{rG_p(s)e^{G_p(r)}}{G_p(r)^{-p+1}},
\]
or equivalently
\[
 r\frac{e^{G_p(2r-s)}}{G_p(2r-s)^{-p}}\left[\frac{2r-s}{r}- \frac{G_p(2r-s)^{-p+1}e^{-G_p(2r-s)}}{G_p(r)^{-p+1}e^{-G_p(r)}}\right]\geq r\frac{e^{G_p(s)}}{G_p(s)^{-p}}\left[-\frac{s}{r}+\frac{G_p(s)^{-p+1}e^{-G_p(s)}}{G_p(r)^{-p+1}e^{-G_p(r)}}\right].
\]
Recall that by item ii) of Lemma \ref{lemma properties of function}, we have that $\frac{e^{G_p(2r-s)}}{G_p(2r-s)^{-p}}\geq \frac{e^{G_p(s)}}{G_p(s)^{-p}}$, since $2r-s\geq s$. Moreover, item iii) of Lemma \ref{lemma properties of function} gives that
\[
G_p(r)^{-p+1}e^{-G_p(r)}\geq \frac{1}{2}G_p(s)^{-p+1}e^{-G_p(s)}+\frac{1}{2}G_p(2r-s)^{-p+1}e^{-G_p(2r-s)},
\]
which after rearranging is equivalent to
\[
\frac{2r-s}{r}- \frac{G_p(2r-s)^{p+1}e^{-G_p(2r-s)}}{G_p(r)^{-p+1}e^{-G_p(r)}}\geq -\frac{s}{r}+\frac{G_p(s)^{-p+1}e^{-G_p(s)}}{G_p(r)^{-p+1}e^{-G_p(r)}}.
\]

Finally, it is easy to see that the inequality becomes equality when $f=c\chi_{[0,d]}$, which indeed corresponds to the $x=0$ (limiting) case.
\end{proof}

\begin{remark}\label{remark after main theorem}
    Since the function $f$ that we consider is decreasing we have -- for any positive Borel measure $\mu$ -- that
    \[
    \frac{\int_0^hfd\mu}{\mu([0,h])}\leq f(0).
    \]
Moreover, by the fact that $N$ is increasing, and denoting by $\Phi$ the distribution function of $\mu$, it is straightforward to check that the function
\[
\frac{1}{s} \int_0^{\Phi^{-1}(s)}N\ d\mu , \text{ for } s>0,
\]
is also increasing. Thus
\[
\int_0^{+\infty}N(t)f(t)d\mu(t) \geq \frac{u}{\Phi(h)} \int_0^{\Phi^{-1}\left(\frac{V}{u}\Phi(h)\right)}N(t)d\mu_{}(t)\geq f(0) \int_0^{\Phi^{-1}\left(\frac{\int_0^{+\infty}f\ d\mu}{f(0)}\right)}N(t)d\mu(t).
\]
  In particular, Theorem \ref{main thm weighted} recovers the classical bound that follows from the reasoning of \cite{fradelizi}, when restricted to even log-concave functions and taking $\mu$ to be the Lebesgue measure, while it is new for the rest measures included in Theorem \ref{main thm weighted}.    
\end{remark}

Choosing $p=\lambda=0$ and $N_q(t)=t^q$, for $q\geq 1$, the expression \eqref{equation main theorem} takes a simple form and gives the following corollary: 

\begin{corollary} \label{corollary lebesgue powers}
    Let $q\geq1$ and $f :[0,+\infty) \to [0,+\infty)$ be an integrable decreasing log-concave function. If $h>0$ is such that 
    $$
    \int_0^h f \leq \frac{1}{2}\int_0^{+\infty} f,
    $$
    then
\begin{equation}
\int_0^{+\infty} t^qf(t)dt\left(\int_0^hf(t)dt\right)^q\geq \frac{1}{q+1}h^q\left(\int_0^{+\infty} f(t)dt\right)^{q+1}.
\end{equation}
    
\end{corollary}

\begin{remark}
    We note that for each choice of convex increasing function $N$, there is a room for improvement of the acceptable ratio $\frac{u}{V}$ for which the inequality holds. Finding the sharp value $\delta_N$ for each such $N$ is certainly of interest and could be important in applications like the extremal slabs of the cube in \cite{BartheKoldobskySlabs}. It was observed in particular that the optimal value for $N_q(t)=t^q$ was an increasing function of $q$.
\end{remark}

\begin{remark}\label{s-concave_lower_bound}

Since the function $c\chi_{_{[0,d]}}$ is also $s$-concave for any $s>0$, all the above inequalities (concerning the lower bound) remain valid in the setting of $s$-concave functions.
\end{remark}

\section{Applications} \label{section applications}

\subsection{Quantitative estimates for norms and entropy}\label{subsection p-norms entropy} Before applying directly Theorem \ref{main thm weighted}, we start by showing how the reduction lemmas --an intermediate step in the previous section-- can be used in other situations. We illustrate this by considering $p$-norms. First, with the notation of Section \ref{subsection reductions}, for any decreasing $M:[0,+\infty)\to [0,+\infty)$ we write
\[
\sup_{f\in \mathcal{F}_{h,u,V}} \int_0^{+\infty} M(t)f(t) dt=VM(0)-\inf_{f\in \mathcal{F}_{h,u,V}} \int_0^{+\infty} \left(M(0)-M(t)\right)f(t) dt.
\]
Thus by Lemma \ref{lemma reduction lower},
\[
\sup_{f\in \mathcal{F}_{h,u,V}} \int_0^{+\infty} M(t)f(t) dt\]
\[=\sup\left\{\int_0^{+\infty} M(t)f(t) dt :  f\in \mathcal{F}_{h,u,V}(\mu)   \text{ such that } f(t)=c \chi_{_{[0,d]}}e^{-at}, \text{  for  } a,c,d\geq 0, d\geq h\right\}.
\]
Let us henceforth denote for brevity \[\mathcal{F}'_{h,u,V}(\mu):=\left\{ f\in \mathcal{F}_{h,u,V}(\mu)   \text{ such that } f(t)=c \chi_{_{[0,d]}}e^{-at}, \text{  for  } a,c,d\geq 0, d\geq h\right\}.\]

Let $f:[0,+\infty)\to [0,+\infty)$ be a decreasing log-concave function and $p,q>1$ be dual exponents. It is known that
\begin{equation}
    \|f\|_{L^p(\mu)}=\sup\left\{ \int_0^{+\infty} f(t) g(t)\ d\mu(t) :\: \|g\|_{L^q(\mu)}=1  \right\}.
\end{equation}
Moreover, the supremum is achieved for $g=\|f\|_p^{1-p}f^{p-1}$, which is a decreasing function. Hence, we can rewrite
\begin{equation}
    \|f\|_{L^p(\mu)}=\sup\left\{ \int_0^{+\infty} f(t) g(t)\ d\mu(t) :\: \|g\|_{L^q(\mu)}=1 , \: \text{decreasing} \right\}.
\end{equation}
 In particular, combining these two facts and interchanging the two suprema, we may write:
 \begin{align*}
      \sup_{f\in \mathcal{F}_{h,u,V}}\|f\|_{L^p(\mu)}&=\sup_{f\in \mathcal{F}_{h,u,V}}\sup_{\|g\|_{L^q(\mu)}=1 , \: \text{decreasing}} \int_0^{+\infty} f(t) g(t)\ d\mu(t)\\
      &=\sup_{\|g\|_{L^q(\mu)}=1 , \: \text{decreasing}}\sup_{f\in \mathcal{F}_{h,u,V}} \int_0^{+\infty} f(t) g(t)\ d\mu(t)\\
      &= \sup_{\|g\|_{L^q(\mu)}=1 , \: \text{decreasing}}\sup_{f\in \mathcal{F}'_{h,u,V}} \int_0^{+\infty} f(t) g(t)\ d\mu(t)\\
      &= \sup_{f\in \mathcal{F}'_{h,u,V}}\|f\|_{L^p(\mu)}.
 \end{align*}
Now, setting $\mu$ to be  the Lebesgue measure, we prove Proposition \ref{p:p-norm}.
\begin{proof}[Proof of Propositon \ref{p:p-norm}]
    The left-hand side of \eqref{p-norms estimate} equals $\|f\|_p^p$. according to the discussion above, an upper bound on this quantity --for any decreasing log-concave function-- will follow from the following optimization problem:
    \[
    \sup\left\{ \int_0^{+\infty} f(t)^p\ dt: f\in \mathcal{F}'_{h,u,V}\right\}.
    \]
    Note also that the right-hand of \eqref{p-norms estimate} equals $\frac{u^{p-1}}{h^{p-1}}V$. Let $f(t)=ce^{-at}\chi_{_{[0,d]}}(t)\in \mathcal{F}'_{h,u,V}$. Then (see also the proof of Lemma \ref{lemma optimization})
    \[
    u=\frac{c}{a}\left(1-e^{-ha}\right), V=\frac{c}{a}\left(1-e^{-da}\right) \text{ and } \|f\|_p^p=\frac{c^p}{pa}\left(1-e^{-pad}\right) .
    \]
    As in Lemma \ref{lemma optimization} (note that in this case it coincides with the parametrization in \cite[Lemma 5]{BartheKoldobskySlabs}) we set 
    \[
    x:=\frac{(1-e^{-ah})V}{u}\in [0,1],
    \]
    and thus $d=\frac{1}{a}(-\log(1-x))$, $a=\frac{1}{h}(-\log(1-\frac{u}{V} x))$ and $c = \frac{aV}{x}$. In particular
    \[
    \|f\|_p^p = \frac{V^p}{h^{p-1}p}\frac{1}{x^p}\left(1 - (1-x)^p\right)\left(-\log\left(1-\frac{u}{V} x\right)\right)^{p-1}.
    \]
    Set $\theta:=\frac{u}{V}$ and consider the function $H_{p}(x) := \frac{1}{x^p}\left(1 - (1-x)^p\right)\bigl(-\log(1-\theta x)\bigr)^{p-1}$. We will show that there exists $\theta_p$, such that for all $\theta\leq \theta_p$, then $H_p(x)\leq p\theta^{p-1}=\lim_{x\to 0^+}H_p(x)$.
    
    Let $k(y):=-\frac{\log(1-y)}{y}$,  with $y\in (0,1)$. Observe that it is increasing and that the desired inequality becomes:
\[
\frac{1 - (1-x)^p}{x}k(\theta x)^{p-1}\leq p \iff \theta \leq \frac{1}{x} k^{-1}\left(\left[\frac{px}{1 - (1-x)^p}\right]^{\frac{1}{p-1}}\right).
\]
The claim now holds for \begin{equation}\label{theta_p in p-norms}\theta_p:=\inf_{x\in (0,1]}\frac{1}{x} k^{-1}\left(\left[\frac{px}{1 - (1-x)^p}\right]^{\frac{1}{p-1}}\right).\end{equation}

Since on $(0,1]$ the expression inside the infimum is continuous and positive, it is only left to check the limit at $0$, in order to ensure that $\theta_p>0$. Indeed we have

\[
\lim_{x\to 0^+} \frac{1}{x} k^{-1}\left(\left[\frac{px}{1 - (1-x)^p}\right]^{\frac{1}{p-1}}\right)=\lim_{x\to 0^+} \frac{1}{p-1}\frac{p^{\frac{1}{p-1}}\left[\frac{x}{1 - (1-x)^p}\right]^{\frac{2-p}{p-1}}\left[\frac{1}{1 - (1-x)^p}-\frac{px(1-x)^{p-1}}{(1 - (1-x)^p)^2}\right]}{k'\left(k^{-1}\left(\left[\frac{px}{1 - (1-x)^p}\right]^{\frac{1}{p-1}}\right)\right)}
\]
\[
=\frac{2p^{\frac{1}{p-1}}}{p-1}\lim_{x\to 0^+}\frac{1}{\left[\frac{1 - (1-x)^p}{x}\right]^{\frac{p}{p-1}}}\times \frac{1 - (1-x)^p-px(1-x)^{p-1}}{x^2}=1,
\]
which concludes the proof. The equality cases are evident and they correspond to $a=0$ or, equivalently, $x=0$.
\end{proof}

It is not hard to see that $\theta_p\geq \frac{1}{2}$ when $1<p\leq 2$. In particular, there is a universal lower bound for $\theta_p$ as $p\to 1^+$, while the inequality becomes an equality for $p=1$. Therefore, differentiating at $p=1$ we get the following:
\begin{corollary}
     Let $f:[0,+\infty)\to [0,+\infty)$ be a decreasing log-concave function. There exists $\theta\geq \frac{1}{2}$, such that if $h>0$ satisfies
     \[
     \int_0^hf(t)dt\leq \theta \int_0^{+\infty} f(t)dt, 
    \]
    then
\begin{equation}
    \int_0^{+\infty} f(t)\log(f(t)) dt \leq \log \left(\frac{1}{h}\int_0^hf(t)dt\right)\int_0^{+\infty} f(t)dt.
\end{equation}

\end{corollary}

The bounds proved in this subsection can be thought of as improvements of the classical Renyi and Shannon entropy bounds for decreasing log-concave probability densities. 

\subsection{Probabilistic interpretations} \label{subsection probabilistic}
In this part of the section we would like to emphasize on the probabilistic nature of the main results. An example of probabilistic application of such results can be found in the recent \cite{Tkocz}, where the main result of \cite{BartheKoldobskySlabs} was used.

\subsubsection{Log-concave random variables} Let $X$ be an even log-concave random variable. That is, a random variable having an even and log-concave density $f$. Note that $f$ is decreasing on $[0,+\infty)$. Hence, Theorem \ref{main theorem intro} can be applied to this function for all $h\leq  \mathrm{median}(|X|)=:m_{|X|}$. This gives in particular
\[
\frac{1}{2}\mathbb{E}N(X)\geq \frac{1}{2m_{|X|}}\int_0^{2m_{|X|}}N(t)\ dt,
\]
for all $N:\R \to \R$ even and convex. Moreover, Theorem \ref{t:upper_bound intro} implies that 
\[
\frac{1}{2}\mathbb{E}X^2\leq \frac{2m_{|X|}^2}{(\log2)^2}.
\]
We remark that not the whole strength of Theorem \ref{t:upper_bound intro} is used. See also Section \ref{anticoncentration section}.

Let now $N:[0,+\infty)\to [0,+\infty)$ convex and increasing. Denote $\tilde{N}(t):=\frac{1}{t}\int_0^tN$, which is increasing as well. Following the previous discussion, we obtain the following estimate for the median of a log-concave random variable on $[0,+\infty)$:
\begin{corollary} \label{corollary median}
    Let $X$ be an even log-concave random variable. Then 
    \begin{equation} \label{median estimates}
      \frac{\log(2)}{2}\|X\|_2  \leq m_{|X|} \leq \frac{1}{2} \tilde{N}^{-1}\left( \frac{1}{2}  \mathbb{E}N(X)\right).
    \end{equation}
\end{corollary}

\begin{remark}
    If $X$ is an even $s$-concave random variable, then Theorem \ref{theorem s-conc square} gives the following bound 
    \begin{equation}
        m_{|X|}\geq \sqrt{\frac{(2s+1)(3s+1)}{2}}\frac{1-2^{-\frac{s}{s+1}}}{s}\|X\|_2,
    \end{equation}
    which also recovers the left-hand side of \eqref{median estimates} as $s\to 0^+$.
\end{remark}

We now consider another important notion in probability theory; the Laplace transform of a random variable. Let $X$ be an even log-concave random variable on $\R$. Its Laplace transform $\Lambda_X : \R \to [0,+\infty)$ is given by
\[
\Lambda_X(s)=\mathbb{E}e^{sX}.
\]

Since $X$ is even, it is readily verifiable that $\Lambda_X$ is even as well. Hence, it suffices to estimate $\Lambda_X$ on $[0,+\infty)$. Note also that by evenness
\[
\Lambda_X(s)=\mathbb{E}\frac{e^{sX}+e^{-sX}}{2}.
\]
For any $s>0$ let $N_s(t):=\frac{e^{st}+e^{-st}}{2}$. This function is even, convex and increasing on $[0,+\infty)$. Therefore, applying Theorem \ref{main theorem intro} for the above choice of $N$, the density $f$ and $h=m_{|X|}$, we obtain the following estimate for the Laplace transform:
\begin{corollary}\label{corollary laplace}
Let $X$ be an even log-concave random variable. Then
    \begin{equation}\label{laplace transform estimate}
        \Lambda_X(s)\geq \frac{e^{2sm_{|X|}}-e^{-2sm_{|X|}}}{4sm_{|X|}}\cdot 
    \end{equation}
\end{corollary}

\subsubsection{An improvement of Jensen's inequality} \label{section jensen}
Let $X$ be a random variable on $[0,+\infty)$ with log-concave tails, that is, the function $t \mapsto \mathbb{P}(X\geq t)$ is log-concave. We will see how Theorem \ref{main theorem intro} can be translated to this setting yielding, under mild assumptions, a quantitative improvement of Jensen's inequality.

\begin{proof}[Proof of Corollary \ref{corollary jensen}]
    By the assumptions, the (decreasing) function $f(t)=\mathbb{P}(X\geq t)$ is log-concave. Moreover, since $N'$ is convex and increasing, we can apply Theorem \ref{main theorem intro} for $f$, $N'$ and the proper $h$ to get that
    \begin{equation}
        \int_0^{+\infty}N'(t)f(t)\ dt\geq \frac{\int_0^hf(t)\ dt}{h}\int_0^{\frac{h}{\int_0^hf(t)\ dt}\int_0^{+\infty}f(t)\ dt}N'(t)\ dt.
    \end{equation}
    The claim now follows by integration.
\end{proof}
The fact that Corollary \ref{corollary jensen} improves Jensen's inequality $\mathbb{E}N(X)\geq N\left(\mathbb{E}X \right)$, can be seen from Remark \ref{remark after main theorem}.

\section{Geometric applications}\label{section geometric}
Let $K\subset \R^n$ be a symmetric convex body and $\theta \in \s^{n-1}$. As already discussed, the function $f(t)=|K\cap \{x\in \mathbb{R}^n:\langle x,\theta\rangle=t\}|_{n-1}$ is $\frac{1}{n-1}$-concave. Moreover, by the symmetry of $K$, it is decreasing on $[0,+\infty)$. In particular, the main theorems can be applied for this choice of $f$. Specifically, the main parameters take the following form:
\begin{equation}
    \int_0^{+\infty}f(t)\ dt=\frac{1}{2}|K|,\ \int_0^hf(t)\ dt=\frac{1}{2}|K\cap \{x\in \mathbb{R}^n: |\langle x,\theta\rangle|\leq h\}|,
\end{equation}
and, for any function $N:\mathbb{R}\to \mathbb{R}$ even
\begin{equation}
    \int_0^{+\infty} N(t)f(t)\ dt = \frac{1}{2}\int_{K}N \left( \langle x,\theta\rangle \right)\ dx.
\end{equation}
Let $\mathrm{Slab}_K(\theta,h):=|K\cap \{x\in \mathbb{R}^n: |\langle x,\theta\rangle|\leq h\}|$ denote the volume of the slab of width $2h$, orthogonal to $\theta$. Then, Theorem \ref{main theorem intro} takes the following form, which can be thought of as quantitative improvement of \cite[Theorem 1]{fradelizi} in the symmetric case, where the sections are now replaced by the slabs of the convex body (for the symmetric case and sections see also \cite{ball}).

\begin{corollary} \label{cor convex bodies lower}
    Let $N:[0,+\infty)\to [0,+\infty)$ be increasing and convex. Let $K\subset \R^n$ be a symmetric convex body with $|K|=1$ and $\theta \in \s^{n-1}$. If $h>0$ is such that
    \begin{equation} 
        \mathrm{Slab}_K(\theta,h)\leq \frac{1}{2},
    \end{equation}
    then
    \begin{equation}
        \int_{K}N \left( \langle x,\theta\rangle \right)\ dx\geq \frac{\mathrm{Slab}_K(\theta,h)}{h}\int_0^{\frac{h}{\mathrm{Slab}_K(\theta,h)}}N(t) \ dt.
    \end{equation}
\end{corollary}
Moreover, Theorem \ref{t: Barthe Koldobsky functional} of Barthe--Koldobsky in conjunction with Theorem \ref{theorem upper 1/n intro} gives the following. Observe that
$$
1-\left(\frac{n}{n+2}\right)^n \geq \frac{3}{4}
$$
for every $n\geq 2$.    
\begin{corollary}\label{cor convex bodies upper}
    Let $K\subset \R^n$ be a symmetric convex body with $|K|=1$ and $\theta \in \s^{n-1}$. If $h>0$ is such that
    \begin{equation} 
        \mathrm{Slab}_K(\theta,h)\leq \frac{3}{4},
    \end{equation}
    then
    \begin{equation}
     \frac{1}{3}\left(\frac{h}{\mathrm{Slab}_K(\theta,h)}\right)^2   \leq \int_K \langle x,\theta\rangle^2\ dx\leq \frac{2h^2}{ \left(n+1\right)\left(n+2\right)\left(1-\left(1-\mathrm{Slab}_K(\theta,h)\right)^{\frac{1}{n}}\right)^2} \cdot 
    \end{equation}
\end{corollary}

A convex body $K\subset \mathbb{R}^n$ is called isotropic, if it has volume one, center of mass at the origin, i.e.,~$\int_Kx\ dx=0$ and there exists $L_K>0$ such that:
\begin{equation} \label{isotropicity}
    \int_K\langle x, \theta\rangle^2\ dx=L_K^2, \text{ for each } \theta\in \s^{n-1}.
\end{equation}

If $K\subset \mathbb{R}^n$ is a convex body, there always exists $T\in GL_n$ so that $TK$ is isotropic (see e.g. \cite[Proposition~2.3.3]{convexbook}. It is then said that $T$ puts $K$ in isotropic position. A central problem in convex geometry for decades was to establish an upper bound for $L_K$ independent of the dimension. This was resolved very recently in \cite{guan,klartaglehec} (see also \cite{bizeul}). An equivalent formulation, in terms of a lower bound on $|K\cap H|$, for any hyperplane $H$, comes from the results of Hensley (see also \cite[Corollary 3.2]{MilPa89}) which say that $L_K$ is equivalent up to universal constants to $\frac{1}{|K\cap H|}$. In the symmetric case, Corollary \ref{cor convex bodies upper} gives the following generalization, where we use symmetric slabs rather than sections:
\begin{corollary} \label{corollary isotropic}
    Let $K\subset \R^n$ be a symmetric isotropic convex body and $\theta \in \s^{n-1}$. If $h>0$ is such that
    \begin{equation} 
        \mathrm{Slab}_K(\theta,h)\leq \frac{3}{4},
    \end{equation}
    then
    \begin{equation}
        \frac{1}{2\sqrt{3}} \frac{2h}{\mathrm{Slab}_K(\theta,h)}\leq L_K \leq \frac{1}{\sqrt{2}} \frac{2h}{ \sqrt{\left(n+1\right)\left(n+2\right)}\left(1-\left(1-\mathrm{Slab}_K(\theta,h)\right)^{\frac{1}{n}}\right)}.
    \end{equation}
\end{corollary}

\subsection{Floating bodies}
Let $K\subset\mathbb{R}^n$ be a convex body and $\delta \in (0,1)$. The convex floating body, introduced in \cite{SW90}, is defined as
\begin{equation}
K_\delta = \cap_{\theta\in \s^{n-1}} \{x \in K : \esc{x,\theta} \leq m_\delta (\theta)\},
\end{equation}
where $m_\delta : \s^{n-1} \to \R$ is such that
\begin{equation}\label{e:floating}
|\{x \in K : \esc{x,\theta} \leq m_\delta (\theta)\}| = (1-\delta)|K|.
\end{equation}
Assuming that $K$ is symmetric it is easy to see that $m_\delta(\theta)=m_\delta(-\theta)$ and in particular
\begin{equation}
    K_\delta = \cap_{\theta\in \s^{n-1}} \{x \in K : |\esc{x,\theta}| \leq m_\delta (\theta)\}.
\end{equation}

It was proved by Fresen \cite{Fresen12} that for any isotropic convex body $K\subset\R^n$ and $\delta \in (0,e^{-1})$,
\begin{equation}\label{fresen inclusion}
(e^{-1}- \delta) L_KB_2^n \subset K_\delta \subset 17\log (\delta^{-1})L_KB_2^n.
\end{equation}
The following is immediate from \eqref{fresen inclusion}, after putting $K$ in isotropic position, and using the fact that $|B_2^n|^{\frac{1}{n}}$ is \emph{equivalent} to $\frac{1}{\sqrt{n}}$ up to an absolute constant.
\begin{corollary}[Fresen]
    Let $K\subset\mathbb{R}^n$ be an isotropic convex body and $\delta\in (0,e^{-1})$. Then,
    \begin{equation}\label{fresen quotient}
        c_1\frac{e^{-1}-\delta}{\sqrt{n}}L_K\leq |K_\delta|^{\frac{1}{n}} \leq c_2\frac{\log(\frac{1}{\delta})}{\sqrt{n}}L_K, 
    \end{equation}
    where $c_1$ and $c_2$ are universal constants.
\end{corollary}

Let $\delta_n = 1-\left(\frac{n}{n+2}\right)^n$.  In a similar vein, if $K$ is isotropic and symmetric, for $\frac{1}{2}\geq \delta\geq \frac{1-\delta_n}{2}$, Corollary \ref{corollary isotropic} gives that for each $\theta\in \s^{n-1}$
\begin{equation}\label{bound for m}
    L_K\sqrt{\frac{(n+1)(n+2)}{2}}\left( 1-2^{\frac{1}{n}}\delta^{\frac{1}{n}} \right)\leq m_\delta(\theta)\leq \sqrt{3}L_K(1-2\delta),
\end{equation}
since for any $\delta\leq \frac{1}{2}$ and $\theta\in \s^{n-1}$, 
$$
\mathrm{Slab}_K(\theta,m_\delta(\theta))=\left(1-2\delta\right).
$$
In particular, the following holds complementing \eqref{fresen inclusion}.
 \begin{lemma}
     Let $K\subset\mathbb{R}^n$ be a symmetric, isotropic convex body and $\delta\in [\frac{1-\delta_n}{2},\frac{1}{2}]$. Then
\begin{equation}\label{inclusion of floatin}
    L_K\sqrt{\frac{(n+1)(n+2)}{2}}\left( 1-2^{\frac{1}{n}}\delta^{\frac{1}{n}} \right)B_2^n\subset K_\delta \subset \sqrt{3}L_K(1-2\delta)B_2^n .
\end{equation}
\end{lemma}

One immediately deduces the following:
\begin{corollary}
    Let $K\subset\mathbb{R}^n$ be a symmetric, isotropic convex body and $\delta\in [\frac{1-\delta_n}{2},\frac{1}{2}]$. Then
    \begin{equation}\label{quotient floatin}
    c'_1\sqrt{\frac{(n+1)(n+2)}{2n}}\left( 1-2^{\frac{1}{n}}\delta^{\frac{1}{n}} \right)L_K\leq |K_\delta|^{\frac{1}{n}} \leq c'_2\frac{\sqrt{3}(1-2\delta)}{\sqrt{n}}L_K, 
    \end{equation}
    where $c'_1$ and $c'_2$ are universal constants.
\end{corollary}

\section{Final Remarks} \label{section discussion}
\subsection{An Anticoncentration type phenomenon} \label{anticoncentration section} In probability theory, anticoncentration bounds are \textit{upper} bounds on the probability that a random variable is close to $0$ or equivalently, \textit{lower} bounds guaranteeing that it stays away from $0$ with high probability. A well known instance of such bound is the Paley-Zygmund inequality stating that for an a.e. non-negative random variable $X$ of finite second moment and $h\in [0,1]$, 
$$
\mathbb{P}(X>h\mathbb{E}X)\geq (1-h)^2\frac{\left(\mathbb{E}X\right)^2}{\mathbb{E}X^2}.
$$ 

We would like to briefly stress the ``anticoncentration" nature of Theorems \ref{t:upper_bound intro} and \ref{theorem upper 1/n intro}. Staying in the $\log$-concave setting for simplicity, assume that $X$ is a random variable with decreasing log-concave density $f$ and $X\geq 0$ a.e. Then, Theorem \ref{t:upper_bound intro} states the following (see Section \ref{subsection probabilistic})
\begin{equation} \label{eq anticoncentration}
\mathbb{P}(X\leq h)\leq 1-e^{-\sqrt{3}} \implies \mathbb{P}(X\leq h) \leq 1-e^{-\sqrt{2}\frac{h}{\|X\|_2}} .
\end{equation}

In order for this to provide something meaningful, it is still left to confirm that $1-e^{-\sqrt{3}} \geq 1-e^{-\sqrt{2}\frac{h}{\|X\|_2}}$, which is not possible in such generality.

\subsection{Towards a general theorem} Let $\psi$ be an increasing function defined on $[0,+\infty)$, $\mu$ a positive Borel measure on $[0,+\infty)$ and consider the class
\begin{equation}\label{general class}
    \mathcal{F}_{h,u,V}^\psi(\mu):= \{f:[0,+\infty)\to [0,+\infty), \psi\text{-concave, decreasing}, \int_0^hf\ d\mu=u, \int_0^{+\infty}f \ d\mu=V  \}.
\end{equation}
As it may be evident from Lemma \ref{lemma reduction upper s} (but also Lemmas \ref{lemma reduction lower} and \ref{lemma reduction upper}), optimizing functionals of the form
\[
\int N f \ d\mu,
\]
over the family $\mathcal{F}_{h,u,V}^\psi(\mu)$ (for a monotone function $N$) can be reduced to the much smaller families of either $\psi$-linear or the combination of constants with $\psi$-linear. We are thus naturally led to the following:

\begin{question}  Can Theorem \ref{main thm weighted} and the corresponding upper bound be proven for general $\mu, \psi$ and $N$?\end{question}

In view of the applications discussed in Sections \ref{section applications} and \ref{section geometric}, many sub-cases of our question seem already interesting. For instance, obtaining a positive answer in the case where $N(t) = t^2$ and $\mu$ is the Lebesgue measure would yield a new result, even in the case of $h\to 0^+$. It is also natural to wonder about the optimal upper bound when $\psi=\log$ and $\mu$ is the Lebesgue measure, thereby generalizing Theorem \ref{t:upper_bound intro}. Furthermore, extending Theorem \ref{main thm weighted} for more general measures remains an open direction. An instance not addressed in this work is when $\mu$ is a discrete measure, which might be of independent interest; to our knowledge, results of this nature are not available in that setting. 

Finally, finding sharp constants of the acceptable ratio $\frac{u}{V}$ (in the terminology of Section \ref{s:optimization}) is an important first step toward understanding the behavior of the inequality in the rest regimes of $\frac{u}{V}$. Apart from novel inequalities and deeper understanding of the behavior of one-dimensional log-concave functions, such result can bring new insights to the geometric question of finding which slabs --of intermediate volume-- are extremal.


\bibliographystyle{alpha}
\bibliography{slabs}

\section*{Appendix}\label{appendix}

\begin{proof}[Proof of Theorem \ref{theorem s-conc square}]
The proof will be along the lines of the proof of Theorem \ref{t:upper_bound intro}, employing now Lemma \ref{lemma reduction upper s}. Let $f\in \mathcal{F}^s_{h,u,V}$. We rewrite $f$ as
\[
f(t)=c\chi_{[0,d]}(t)+c\left(1-\lambda(t-d)\right)^{1/s}\chi_{[d,d+\frac{1}{\lambda}]}(t),
\]
where $d\leq h\leq d+\frac{1}{\lambda}$. By the assumptions on $f$ it holds that
\[
u=cd+\frac{c}{\lambda\left(\frac{1}{s}+1\right)}-\frac{c}{\lambda\left(\frac{1}{s}+1\right)}\left[1-\lambda(h-d)\right]^{\frac{1}{s}+1} \quad \text{and}\quad
V=cd+\frac{c}{\lambda\left(\frac{1}{s}+1\right)}.
\]
Setting $x=d\lambda$, we rewrite the parameters as follows:
\[
c=\frac{V}{d+\frac{1}{\lambda\left(\frac{1}{s}+1\right)}}\iff cd=\frac{V}{1+\frac{1}{x\left(\frac{1}{s}+1\right)}}\iff \frac{c}{\lambda \left(\frac{1}{s}+1\right)}=\frac{V}{x\left(\frac{1}{s}+1\right)+1}.
\]
So 
\[
u=V-\frac{V}{x\left(\frac{1}{s}+1\right)+1}\left[1-\lambda(h-d)\right]^{\frac{1}{s}+1}.
\]
Denoting $\Delta:= 1-\frac{u}{V}\in (0,1)$, we get that
\[
\Delta\left(x\left(\frac{1}{s}+1\right)+1\right)=\left[1-\lambda(h-d)\right]^{\frac{1}{s}+1}\iff \lambda =\frac{x+1-\left[\Delta\left(x\left(\frac{1}{s}+1\right)+1\right)\right]^{\frac{s}{s+1}}}{h},
\]
and
\[
d=\frac{hx}{x+1-\left[\Delta\left(x\left(\frac{1}{s}+1\right)+1\right)\right]^{\frac{s}{s+1}}} \cdot
\]
For any $p>0$, if $N(t)=t^p$ we have
\begin{align*}
    \int_0^{+\infty}N(t)f(t)\ dt&= \int_0^{+\infty}t^pf(t)dt\\
    &=c\int_0^dt^pdt+c\int_d^{d+\frac{1}{\lambda}}t^p\left(1-\lambda(t-d)\right)^{1/s}dt\\
    &=\frac{cd^{p+1}}{p+1}+\frac{c}{\lambda^{p+1}}\int_0^1\left(1-t+d\lambda\right)^{p}t^{\frac{1}{s}}dt.
\end{align*}
In particular, for $p=2$
\begin{equation}\label{e:appendix1}
    \begin{split}
    \int_0^{+\infty}t^2f(t)\ dt&=  \frac{cd^{3}}{3}+\frac{c}{\lambda^{3}}\int_0^1\left(1-t+d\lambda\right)^{2}t^{\frac{1}{s}}dt\\
    &=\frac{1}{3}\frac{V}{1+\frac{1}{x\left(\frac{1}{s}+1\right)}}\left(\frac{hx}{x+1-\left[\Delta\left(x\left(\frac{1}{s}+1\right)+1\right)\right]^{\frac{s}{s+1}}}\right)^2\\
    &\quad+\frac{V\left(\frac{1}{s}+1\right)}{x\left(\frac{1}{s}+1\right)+1}\left(\frac{h}{x+1-\left[\Delta\left(x\left(\frac{1}{s}+1\right)+1\right)\right]^{\frac{s}{s+1}}}\right)^2\int_0^1\left(1-t+x\right)^{2}t^{\frac{1}{s}}dt.
    \end{split}
\end{equation}
The expression \eqref{e:appendix1} equals at $x=0$ 
\[
\frac{2Vh^2}{\left(1-\Delta^{\frac{s}{s+1}}\right)^2\left(\frac{1}{s}+2\right)\left(\frac{1}{s}+3\right)}.
\]
Thus, it suffices to show that for $\Delta$ sufficiently large that
\begin{align*}
    \frac{1}{3}&\frac{V}{1+\frac{1}{x\left(\frac{1}{s}+1\right)}}\left(\frac{hx}{x+1-\left[\Delta\left(x\left(\frac{1}{s}+1\right)+1\right)\right]^{\frac{s}{s+1}}}\right)^2\\
    &+\frac{V\left(\frac{1}{s}+1\right)}{x\left(\frac{1}{s}+1\right)+1}\left(\frac{h}{x+1-\left[\Delta\left(x\left(\frac{1}{s}+1\right)+1\right)\right]^{\frac{s}{s+1}}}\right)^2\int_0^1\left(1-t+x\right)^{2}t^{\frac{1}{s}}dt\\
    &\leq \frac{2Vh^2}{\left(1-\Delta^{\frac{s}{s+1}}\right)^2\left(\frac{1}{s}+2\right)\left(\frac{1}{s}+3\right)}
\end{align*}
or equivalently,
\begin{align*}
    \frac{1}{3}&\frac{1}{1+\frac{1}{x\left(\frac{1}{s}+1\right)}}\left(\frac{x}{x+1-\left[\Delta\left(x\left(\frac{1}{s}+1\right)+1\right)\right]^{\frac{s}{s+1}}}\right)^2\\
    &+\frac{\left(\frac{1}{s}+1\right)}{x\left(\frac{1}{s}+1\right)+1}\left(\frac{1}{x+1-\left[\Delta\left(x\left(\frac{1}{s}+1\right)+1\right)\right]^{\frac{s}{s+1}}}\right)^2\int_0^1\left(1-t+x\right)^{2}t^{\frac{1}{s}}dt\\
    &\leq \frac{2}{\left(1-\Delta^{\frac{s}{s+1}}\right)^2\left(\frac{1}{s}+2\right)\left(\frac{1}{s}+3\right)}.
\end{align*}

Considering the function
\begin{align*}
    H(x):=&\left(\frac{x^3}{3}+\int_0^1\left(1-t+x\right)^{2}t^{\frac{1}{s}}dt\right)\left(1-\Delta^{\frac{s}{s+1}}\right)^2\left(\frac{1}{s}+1\right)\left(\frac{1}{s}+2\right)\left(\frac{1}{s}+3\right)\\
&-\left(x+1-\left[\Delta\left(x\left(\frac{1}{s}+1\right)+1\right)\right]^{\frac{s}{s+1}}\right)^22\left(x\left(\frac{1}{s}+1\right)+1\right),
\end{align*} 
which is $0$ at $x=0$ (as we observed above), it is enough to prove that $H'(x)\leq 0$. Moreover, we have that
\begin{align*}
    H'(x)=&\left(x^2+2\int_0^1\left(1-t+x\right)t^{\frac{1}{s}}dt\right)\left(1-\Delta^{\frac{s}{s+1}}\right)^2\left(\frac{1}{s}+1\right)\left(\frac{1}{s}+2\right)\left(\frac{1}{s}+3\right)\\
&-\left(x+1-\left[\Delta\left(x\left(\frac{1}{s}+1\right)+1\right)\right]^{\frac{s}{s+1}}\right)^22\left(\frac{1}{s}+1\right)\\
&-4\left(x\left(\frac{1}{s}+1\right)+1\right)\left(x+1-\left[\Delta\left(x\left(\frac{1}{s}+1\right)+1\right)\right]^{\frac{s}{s+1}}\right)\left(1-\Delta^{\frac{s}{s+1}}\left(x\left(\frac{1}{s}+1\right)+1\right)^{-\frac{1}{s+1}}\right)
\end{align*}
Note that
\begin{align*}
    H'(0)&=\left(2\int_0^1\left(1-t\right)t^{\frac{1}{s}}dt\right)\left(1-\Delta^{\frac{s}{s+1}}\right)^2\left(\frac{1}{s}+1\right)\left(\frac{1}{s}+2\right)\left(\frac{1}{s}+3\right)\\
&\quad- \left(1-\Delta^{\frac{s}{s+1}}\right)^22\left(\frac{1}{s}+1\right)-4\left(1-\Delta^{\frac{s}{s+1}}\right)^2\\
&=0.
\end{align*}
Hence, it suffices to show that $H''(x)\leq 0$ so $H'$ is decreasing. In this regard,
\begin{align*}
    H''(x)=&\left(x+\int_0^1t^{\frac{1}{s}}dt\right)2\left(1-\Delta^{\frac{s}{s+1}}\right)^2\left(\frac{1}{s}+1\right)\left(\frac{1}{s}+2\right)\left(\frac{1}{s}+3\right)\\
&-\left(x+1-\left[\Delta\left(x\left(\frac{1}{s}+1\right)+1\right)\right]^{\frac{s}{s+1}}\right)\left(1-\Delta^{\frac{s}{s+1}}\left(x\left(\frac{1}{s}+1\right)+1\right)^{-\frac{1}{s+1}}\right)4\left(\frac{1}{s}+1\right)\\
&-4\left(\frac{1}{s}+1\right)\left(x+1-\left[\Delta\left(x\left(\frac{1}{s}+1\right)+1\right)\right]^{\frac{s}{s+1}}\right)\left(1-\Delta^{\frac{s}{s+1}}\left(x\left(\frac{1}{s}+1\right)+1\right)^{-\frac{1}{s+1}}\right)\\
&-4\left(x\left(\frac{1}{s}+1\right)+1\right)\left(1-\Delta^{\frac{s}{s+1}}\left(x\left(\frac{1}{s}+1\right)+1\right)^{-\frac{1}{s+1}}\right)^2\\
&-\frac{4\Delta^{\frac{s}{s+1}}}{s}\left(x\left(\frac{1}{s}+1\right)+1\right)\left(x+1-\left[\Delta\left(x\left(\frac{1}{s}+1\right)+1\right)\right]^{\frac{s}{s+1}}\right)\left(x\left(\frac{1}{s}+1\right)+1\right)^{-\frac{1}{s+1}-1}.
\end{align*}
We would like to show that for each $x\geq 0$ it holds $H''(x)\leq 0$. To simplify the notation let $y:=x\left(\frac{1}{s}+1\right)+1\geq 1$ and $q=\Delta^{\frac{s}{s+1}}$. Then the desired inequality is (for all $y\geq 1$)
\begin{align*}
    2&\left(1-q\right)^2\left(\frac{1}{s}+2\right)\left(\frac{1}{s}+3\right)y\\
    &- \left(\frac{s}{s+1}y+\frac{1}{s+1}-qy^{\frac{s}{s+1}}\right)\left(1-qy^{-\frac{1}{s+1}}\right)4\left(\frac{1}{s}+1\right)\\
    &-4\left(\frac{1}{s}+1\right)\left(\frac{s}{s+1}y+\frac{1}{s+1}-qy^{\frac{s}{s+1}}\right)\left(1-qy^{-\frac{1}{s+1}}\right)\\
    &-4y\left(1-qy^{-\frac{1}{s+1}}\right)^2\\
    &-\frac{4q}{s}\left(\frac{s}{s+1}y+\frac{1}{s+1}-qy^{\frac{s}{s+1}}\right)y^{-\frac{1}{s+1}}\leq 0
\end{align*}
Equivalently,
\begin{align*}
    y&\left(1-q\right)^2\left(\frac{1}{s}+2\right)\left(\frac{1}{s}+3\right)\\
    &-4y-\frac{4}{s}+4\frac{s+1}{s}qy^{\frac{s}{s+1}}+4qy^{\frac{s}{s+1}}+\frac{4}{s}qy^{-\frac{1}{s+1}}-4\frac{s+1}{s}q^2y^{\frac{s-1}{s+1}}\\
    &-2y+4qy^{\frac{s}{s+1}}-2q^2y^{\frac{s-1}{s+1}}-\frac{2q}{s+1}y^{\frac{s}{s+1}}-\frac{2q}{s(s+1)}y^{-\frac{1}{s+1}}+\frac{2q^2}{s}y^{\frac{s-1}{s+1}} \leq 0
\end{align*}
or
\begin{align*}
    y&\left(\left(1-q\right)^2\left(\frac{1}{s}+2\right)\left(\frac{1}{s}+3\right)-6  \right)\\
    &-\frac{4}{s}+y^{\frac{s}{s+1}}\left( 4\frac{s+1}{s}q+8q  -\frac{2q}{s+1}    \right)\\
    &+y^{-\frac{1}{s+1}}\left( \frac{4}{s}q-\frac{2q}{s(s+1)} \right)-y^{\frac{s-1}{s+1}}\left(6q^2+\frac{2}{s}q^2  \right)\leq 0.
\end{align*}
In particular, for $y=1$ we would like to have
\begin{align*}
    &\left(\left(1-q\right)^2\left(\frac{1}{s}+2\right)\left(\frac{1}{s}+3\right)-6  \right)-\frac{4}{s}\\
    &+\left( 4\frac{s+1}{s}q+8q  -\frac{2q}{s+1}    \right)\\
    &+\left( \frac{4}{s}q-\frac{2q}{s(s+1)} \right) -\left(6q^2+\frac{2}{s}q^2  \right)\leq 0,
\end{align*}
which is equivalent to
\[
s+1-(4s+2)q+(3s+1)q^2\leq 0\implies q\geq \frac{s+1}{3s+1}.
\]
This is a first restriction for the value of $q$ (or equivalently for $u/V$).

We will assume from now on that this bound holds for $q$. Under this condition, we aim to prove that the following function is decreasing
\begin{align*}
    G(y):=&y^{1+\frac{1}{s+1}} \left(\left(1-q\right)^2\left(\frac{1}{s}+2\right)\left(\frac{1}{s}+3\right)-6  \right)-\frac{4}{s}y^{\frac{1}{s+1}}\\
    &+y\left(\frac{4}{s}q+12q  -\frac{2q}{s+1}    \right)+\left( \frac{4}{s}q-\frac{2q}{s(s+1)} \right)-y^{\frac{s}{s+1}}\left(6q^2+\frac{2}{s}q^2  \right)
\end{align*}
for all large enough $q$, which in turn will prove the theorem. We calculate
\begin{align*}
    G'(y)&=y^{\frac{1}{s+1}}\frac{s+2}{s+1} \left(\left(1-q\right)^2\left(\frac{1}{s}+2\right)\left(\frac{1}{s}+3\right)-6  \right)-\frac{4}{s(s+1)}y^{-\frac{s}{s+1}}\\
    &+\left( \frac{4}{s}q+12q  -\frac{2q}{s+1}    \right)
-y^{-\frac{1}{s+1}}\frac{s}{s+1} \left(6q^2+\frac{2}{s}q^2  \right).
\end{align*}
Note that  $G'(1)\leq 0$ for all $q\geq\frac{s+1}{3s+1}$. Thus, it is enough to see that $G''(y)\leq 0$. Now
\begin{align*}
    G''(y)=&y^{-\frac{s}{s+1}}\frac{s+2}{(s+1)^2} \left(\left(1-q\right)^2\left(\frac{1}{s}+2\right)\left(\frac{1}{s}+3\right)-6  \right)\\
    &+\frac{4}{(s+1)^2}y^{-1-\frac{s}{s+1}}+y^{-1-\frac{1}{s+1}}\frac{s}{(s+1)^2} \left(6q^2+\frac{2}{s}q^2  \right)\leq 0
\end{align*}
if and only if 
\begin{align*}
    \frac{s+2}{(s+1)^2} \left(\left(1-q\right)^2\left(\frac{1}{s}+2\right)\left(\frac{1}{s}+3\right)-6  \right)+\frac{4}{(s+1)^2}y^{-1}+y^{-\frac{2}{s+1}}\frac{s}{(s+1)^2} \left(6q^2+\frac{2}{s}q^2  \right)\leq 0
\end{align*}
Observe that this function is decreasing in $y$, hence it suffices to have the inequality for $y=1$. Putting $y=1$ we get
\begin{align*}
    \frac{s+2}{(s+1)^2} \left(\left(1-q\right)^2\left(\frac{1}{s}+2\right)\left(\frac{1}{s}+3\right)-6  \right)+\frac{4}{(s+1)^2}+\frac{s}{(s+1)^2} \left(6q^2+\frac{2}{s}q^2  \right)\leq 0,
\end{align*}
which is equivalent to
$$
(s+2)\left(1-2q+q^2\right)\left(1+5s+6s^2\right)-6s^3-8s^2+\left(6s^3+2s^2\right)q^2\leq 0.
$$
The latter is a quadratic polynomial in $q$, whose value at $q=1$ is strictly negative. In particular, there exists $q_s<1$ such that the inequality holds for each $q\geq q_s$. Moreover, it is not difficult to see that $\frac{s+1}{3s+1} \geq q_s$ for every $s\geq 0$. Therefore, we have
\begin{equation}
    \frac{u}{V}\leq 1-\left(\frac{s+1}{3s+1}\right)^{\frac{1}{s}+1},
\end{equation}
and we conclude the proof for 
\begin{equation}\label{bound u/V s-conc}
    \delta_s=1-\left(\frac{s+1}{3s+1}\right)^{\frac{1}{s}+1}.
\end{equation}
\end{proof}
\end{document}